\documentclass[12pt]{article}
\usepackage[utf8]{inputenc}
\usepackage[english]{babel}
\usepackage{amsthm}
\usepackage{hyperref}
\usepackage{cleveref}
\usepackage{amsmath, amssymb, mathrsfs}
\usepackage{xypic}
\usepackage{graphicx}
\usepackage{aliascnt}

\newtheorem{theorem}{Theorem}[section]
\newtheorem{corollary}[theorem]{Corollary}
\newtheorem{lemma}[theorem]{Lemma}
\theoremstyle{definition}
\newtheorem{definition}[theorem]{Definition}

\theoremstyle{definition}
\newtheorem*{definition*}{Notational Conventions}

\theoremstyle{remark}
\newtheorem{remark}[theorem]{Remark}

\theoremstyle{definition}
\newtheorem{example}[theorem]{Example}
\begin{document}
\title{The Local Joyal Model Structure}
\author{Nicholas Meadows}
\maketitle

\section*{Introduction}

 The purpose of this paper is develop an analog of the Jardine model structure on simplicial presheaves in which, rather than having the weak equivalences be 'local Kan equivalences', the weak equivalences are 'local Joyal equivalences'. This model structure is called the local Joyal model structure. 

The motivation for the creation of this model structure was to develop a tool for the study of higher-dimensional automata. Higher-dimensional automata are finite cubical sets $K$ that model concurrent processes; each k-cell of a cubical complex represents k processes acting simultaneously, the vertices represent system states, and morphisms in the path category $P(K)(x, y)$ represent execution paths between states x and y. Computing the path category for a cubical set specializes to computing the path category of a simplicial set via the triangulation functor (see \cite{pathcatcomp}). The path category object is naturally understood in terms of the Joyal model structure; Joyal equivalences induce equivalences of path categories, but path categories are not a standard homotopy invariant. The hope is that this model structure will provide a framework to apply homotopy theoretic techniques such as descent theory to get a local to global analysis of the behavior of cellular automata.

Section 1 reviews some facts about quasi-categories and the Joyal model structure that will be used to prove the main theorem of the text. In particular,  the path category and core of quasi-category are described. Joyal equivalences are characterized in a manner compatible with Boolean localization (\Cref{lem1.10}).

Section 2 is devoted to reviewing the technique of Boolean localization, which is essential to proving the existence of the Jardine model structure for simplicial presheaves (as well as the model structure of this paper). Boolean localization states that every Grothendieck topos has a surjective morphism to the topos of sheaves on a complete Boolean algebra. This theorem is proven in \cite{MM}. The article \cite{Boolean} of Jardine gives a proof of the existence of the Jardine model structure based on the technique of Boolean localization. However, Jardine's recent book, \cite{local}, is a more recent and complete exposition of various model structures on simplicial presheaves and their construction.

Section 3 is devoted to proving the existence of the Joyal model structure.  Section 4 describes the corresponding model structure for simplicial sheaves. It also gives a concrete description of the local Joyal model structure on the simplicial sheaves on a Boolean site.

\begin{definition*}
 Write $B(\mathcal{C})$ for the nerve of a small category $\mathcal{C}$. Write $sSet$ for the category of simplicial sets and $Cat$ for the category of small categories. Given simplicial sets $K, Y$, let $hom(K, Y)$ denote the set of morphisms between them. Write $\textbf{hom}(K, Y)$ for the simplicial set whose n-simplices are maps $\Delta^{n} \times K \rightarrow Y$.  Write $P : sSet \rightarrow Cat$ for the left adjoint of $B$; if X is a simplicial set $P(X)$ is called the \textbf{path category} of X. In \cite{Joyal-quasi-cat}, the notation $\tau_{1}$ is used for the left adjoint of $B$ and $\tau_{1}(X)$ is called the fundamental category of $X$. The notation $P$ was chosen due to the theoretical computer science motivations behind this paper (see the introduction to \cite{pathcatcomp}). Write $\pi(X) = P(X)[P(X)]^{-1}$ for the fundamental groupoid of a simplicial set. Finally, $\tau_{0}(K, X)$ will denote \textbf{Joyal's set}, which is defined to be the isomorphism classes in $P(\textbf{hom}(K, X))$.

Denote by $s\textbf{Pre}(\mathscr{C})$ the simplicial presheaves on a Grothendieck site $\mathscr{C}$. Denote by $s\textbf{Sh}(\mathscr{C})$ the simplicial sheaves on $\mathscr{C}$. Write $L^{2} : s\textbf{Pre}(\mathscr{C}) \rightarrow s\textbf{Sh}(\mathscr{C})$ for the sheafification functor (see \cite{Boolean}, \cite{local}). \textbf{Local weak equivalences} are defined to be weak equivalence in the Jardine model structure on simplicial presheaves, as described in \cite[pg. 63-64]{local}. A sectionwise weak equivalence $f: X \rightarrow Y$ of simplicial presheaves on a site $\mathscr{C}$ is a map of simplicial presheaves so that $X(U) \rightarrow Y(U)$ is a weak equivalence for all $U \in \mathscr{C}$, 

\end{definition*}

\section{Preliminaries on Quasi-Categories}
\begin{definition}\label{def1.1}
An \textbf{inner fibration} is a map of simplicial sets which has the right lifting property with respect to all inner horn inclusions $\Lambda_{i}^{n} \subset \Delta^{n}, 0 < i < n$. Say that a simplicial set $X$ is a \textbf{quasi-category} if the map $X \rightarrow *$ is an inner fibration. 
\end{definition}

The existence of the Joyal model structure on simplicial sets is asserted in \cite[Theorem 2.2.5.1]{Lurie} and \cite[Theorem 6.12]{Joyal-quasi-cat}. The fibrant object of this model structure are the quasi-categories.
The cofibrations are monomorphisms. The weak equivalences of this model structure are called \textbf{Joyal equivalences}.  They are defined to be maps $f: A \rightarrow B$, so that for each quasi-category $X$, the map 
$$
\tau_{0}(B, X) \rightarrow \tau_{0}(A, X)
$$
is a bijection. The fibrations of this model structure are called \textbf{quasi-fibrations}. The trivial fibrations are the trivial Kan fibrations. 
\begin{lemma}\label{lem1.2}
The functor $P$ is left adjoint to the nerve functor $B$. Moreover, $P$ preserves finite products.
\end{lemma}
\begin{lemma}\label{lem1.3}
A Joyal equivalence induces an equivalence of path categories.
\end{lemma}
\begin{proof}
 Let $f: X \rightarrow Y$ be a Joyal equivalence of quasi-categories. Then $X \times I $ is a cylinder object for $X$ in the Joyal model structure. Then there exists a map $g: Y \rightarrow X$ and homotopies $f \circ g \sim id_{Y}, g \circ f \sim id_{X}$. Since $P(I) = \pi(\Delta^{1})$ and $P$ preserves finite products, $P(f)$ is an equivalence of categories.  

If $f: X \rightarrow Y$ is a Joyal equivalence, form the diagram
$$
\xymatrix{
X \ar[r]_{i_{X}} \ar[d]^{f} & \mathcal{L}(X) \ar[d]^{\mathcal{L}(f)} \\
Y \ar[r]_{i_{Y}} & \mathcal{L}(Y)}
$$
where the horizontal maps are the natural fibrant replacements for the Joyal model structure (i.e. constructed by taking transfinite composites of pushouts of inner horn inclusions).  By \cite[Lemma 1.6]{Joyal-quasi-cat} the maps $P(\Lambda_{i}^{n}) \rightarrow P(\Delta^{n})$ induced by inner horn inclusions are isomorphism. Thus, since $P$ commutes with products $P(i_{X})$ and $P(i_{Y})$ are isomorphisms. Thus, the first paragraph implies that $P(X) \rightarrow P(Y)$ is an equivalence of categories, as required. 

\end{proof}
\begin{definition}\label{def1.4}
Suppose that X is a quasi-category. Say that 1-simplices $\alpha, \beta: x \rightarrow y$ of X are \textbf{right homotopic}, written by $\alpha \Rightarrow_{R} \beta$, if and only if there exists a 2-simplex with boundary
$$
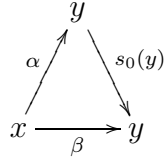
\xymatrix@R=0pt@=10pt@H=6pt{
 & y \ar[ddr]^{s_{0}(y)}& \\
&& \\
x \ar[uur]^{\alpha} \ar[rr]_{\beta} && y \\ }
$$

Similarly, say that $\beta, \alpha$ are \textbf{left homotopic} (written $\beta \Rightarrow_{L} \alpha$) if and only if there exists a 2-simplex with boundary: 
$$
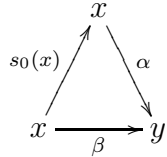
\xymatrix@R=0pt@=10pt@H=6pt{
 & x \ar[ddr]^{\alpha} & \\
&& \\
x \ar[uur]^{s_{0}(x)} \ar[rr]_{\beta} && y \\
}
$$
\end{definition}
\begin{lemma}\label{lem1.5}
If $\alpha$ and $\beta$ are 1-simplices, then $$\alpha \Rightarrow_{R} \beta \; (i)
\iff
\beta \Rightarrow_{R} \alpha \; (ii)
\iff \alpha \Rightarrow_{L} \beta \; (iii)
\iff \beta \Rightarrow_{L} \alpha \; (iv)
$$
If any of the above relations are true then say that $\alpha$ and $\beta$ are \textbf{homotopic}. Moreover, homotopy in this sense is an equivalence relation.
\end{lemma}
\begin{example}\label{exam1.6}
(see \cite[pg. 29-32]{Lurie}) In the case $X$ is a quasi-category, $P(X) = ho(X)$, where $ho(X)$ has following description. It is the category which has objects the vertices of X, and morphisms the homotopy classes of 1-simplices $[\alpha] : x \rightarrow y$ in X. Composition is defined for classes $[\alpha]: x \rightarrow y$, $[\beta] : y \rightarrow z$, by $[d_{1}(\sigma)] = [\beta] \circ [\alpha]$, where $\sigma$ is the 2-simplex depicted in the below diagram
$$ 
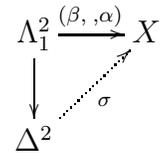
\xymatrix{
\Lambda_{1}^{2} \ar[r]^{(\beta, \;, \alpha)} \ar[d] & X \\
\Delta^{2} \ar@{.>}[ur]_{\sigma} & \\
}
$$
In this category, $s_{0}(x) = id_{x}$.
\end{example}
\begin{lemma}\label{lem1.7}
 Let $I = B\pi(\Delta^{1})$, and $X$ be a quasi-category. A 1-simplex $f$ is invertible in $P(X)$ if and only if there exists a lift in the diagram
$$
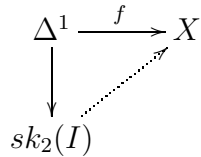
\xymatrix{
\Delta^{1} \ar[r]^{f} \ar[d] & X \\
sk_{2}(I) \ar@{.>}[ur] &  }
$$
\end{lemma}
\begin{proof}
It suffices to produce 2-simplices with boundaries.
$$
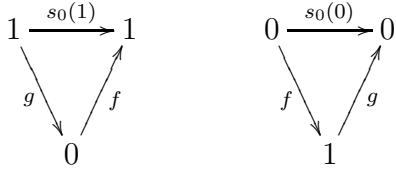
\xymatrix@R=0pt@=10pt@H=6pt{
1 \ar[rr]^{s_{0}(1)} \ar[ddr]_{g} && 1 &&& 0 \ar[rr]^{s_{0}(0)} \ar[ddr]_{f} && 0 \\
&&&&&&& \\
& 0 \ar[uur]_{f} & &&& & 1 \ar[uur]_{g} & \\
}
$$

By symmetry, it suffices to produce a 2-simplex with boundary as depicted on the left. 

Necessity follows from \Cref{exam1.6}. Suppose that $f$ has a left inverse in $P(X)$, $g$, so that $f \circ g$ is right-homotopic to the identity. Consider a map $\sigma : \Lambda_{2}^{3} \rightarrow X$, so that $\sigma_{012}$ expresses $f \circ g$ as a composite of $f, g$ and $\sigma_{123}, \sigma_{023}$ respectively are the 2-simplices with boundaries depicted below: 
$$
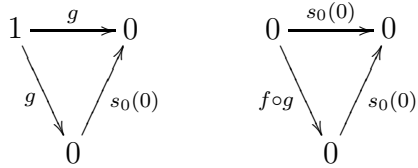
\xymatrix@R=0pt@=10pt@H=6pt{
1 \ar[rr]^{g} \ar[ddr]_{g} && 0 &&& 0 \ar[rr]^{s_{0}(0)} \ar[ddr]_{f \circ g} && 0 \\
&&&&&&& \\
& 0 \ar[uur]_{s_{0}(0)} & &&& & 0 \ar[uur]_{s_{0}(0)} & \\
}
$$
Note that $\sigma_{023}$ expresses the right homotopy between $f \circ g$ and $s_{0}(0)$. Extending $\sigma$ to a 3-simplex $\sigma'$, $d_{2}(\sigma')$ gives the required 2-simplex. 
\end{proof}

The main result of \cite{Joyal1} is as follows: 
\begin{theorem}\label{thm1.8}
A quasi-category $X$ is a Kan complex if and only if its path category $P(X)$ is a groupoid. 
\end{theorem}
\begin{definition}\label{def1.9}
There is a functor $J : Quasi \rightarrow Kan$ where $Quasi$ and $Kan$ are, respectively, the full subcategories of $sSet$ of quasi-categories and Kan complexes. 
\end{definition}
The functoriality of $J$ follows from the fact that the simplices of $J(X)$ are precisely the simplices of $X$ whose edges are invertible in $P(X)$.   
\begin{lemma}\label{lem1.10} A map
 $f: X \rightarrow Y$ is a Joyal equivalence of quasi-categories if and only if
$$
\pi(J\textbf{hom}(K, X)) \rightarrow \pi(J\textbf{hom}(K, Y))
$$
is an equivalence of categories for all finite simplicial sets $K$. 
\end{lemma}
\begin{proof}

 Note that $f$ is a Joyal equivalence if and only if $[K, X] \rightarrow [K, Y]$ is a bijection for each finite simplicial set $K$. By \cite[Poroposition 2.2.5.7]{Lurie}, $[K, X] $ can be identified with $\tau_{0}(K, X)$ which is in bijective correspondence with $\pi_{0}(J\textbf{hom}(K, X))$ by \Cref{thm1.8}. It follows that $f$ is a Joyal equivalence if and only if

$$
\pi_{0}(J\textbf{hom}(K, X)) \rightarrow \pi_{0}(J\textbf{hom}(K, Y))
$$
is a bijection for each finite simplicial set $K$.

By \cite[Proposition 4.26]{Joyal-quasi-cat}, $J$ sends Joyal equivalences of quasi-categories to Joyal equivalences. Thus, since $\textbf{hom}(K, -)$ preserves Joyal equivalences, $J\textbf{hom}(K, X) \rightarrow J\textbf{hom}(K, Y)$ is a Joyal equivalence for all finite simplicial sets $K$. By \Cref{lem1.3}, 
$$
\pi(J\textbf{hom}(K, X)) \rightarrow \pi(J\textbf{hom}(K, Y))
$$
is an equivalence of categories for finite simplicial sets $K$. Combining this with the statement proven in the first paragraph of the proof, the result follows.  
\end{proof}
\section{Preliminaries on Boolean Localization} Given a finite simplicial set $K$ and a simplicial presheaf $X$, write $hom(K, X)$ for the simplicial presheaf $U \mapsto hom(K, X(U))$. Write $\textbf{hom}(K, X)$ for the simplicial presheaf $U \mapsto \textbf{hom}(K, X(U))$. 

\begin{definition}\label{def2.1}
Let $\mathscr{L}, \mathscr{M}$ be Grothendieck topoi. A \textbf{geometric morphism} $p: \mathscr{L} \rightarrow \mathscr{M}$ is a pair of functors $p^{*}: \mathscr{M} \rightarrow \mathscr{L}$, $p_{*} : \mathscr{L} \rightarrow \mathscr{M}$ so that $p^{*}$ preserves finite limits and is left adjoint to $p_{*}$. Call a geometric morphism \textbf{surjective} if and only if it satisfies the following equivalent properties:
\begin{enumerate}
\item{$p^{*}$ is faithful}
\item{$p^{*}$ reflects isomorphisms}
\item{$p^{*}$ reflects monomorphisms}
\item{$p^{*}$ reflects epimorphisms}
\end{enumerate}
\end{definition}
The following theorem is proven in \cite[pg. 515]{MM}, as well as in \cite[Section 3.5]{local}:

\begin{theorem}\label{thm2.2}
 (Barr). Let $\mathscr{L}$ be any Grothendieck topos. Then there exists a surjective geometric morphism $p = (p^{*}, p_{*}): \textbf{Sh}(\mathscr{B}) \rightarrow \mathscr{L}$, so that $\mathscr{B}$ is a complete Boolean algebra. 
\end{theorem}

Such a surjective geometric morphism (from $\textbf{Sh}(\mathscr{B})$) is called a \textbf{Boolean localization} of $\mathscr{L}$.

\begin{definition}\label{def2.3}
Suppose that $i : K \subseteq L$ is an inclusion of finite simplicial sets, and $f: X \rightarrow Y$ is a map of simplicial presheaves. Say that $f$ has the \textbf{local right lifting property} with respect to $i$ if for every commutative diagram

$$
\xymatrix{
 K \ar[r] \ar[d] & X(U)\ar[d]\\
L \ar[r] & Y(U) \\ 
}
$$
there is some covering sieve $R \subseteq hom(-, U), U \in Ob(\mathscr{E})$, such that the lift exists in the diagram
$$
\xymatrix{
K \ar[r] \ar[d] & X(U) \ar[r]^{X(\phi)} & X(V) \ar[d]\\
L \ar[r] \ar@{.>}[urr] & Y(U) \ar[r]_{Y(\phi)} & Y(V)\\
}
$$
for each $\phi \in R$. Similarly, say that $f$ has the sectionwise right lifting property with respect to $i$ if and only if there exists a lifting 
$$
\xymatrix{
 K \ar[r] \ar[d] & X(U)\ar[d]\\
L \ar[r] \ar@{.>}[ur] & Y(U) \\ 
}
$$
for each $U \in Ob(\mathscr{E})$. 
\end{definition}

\begin{definition}\label{def2.4}
Say that a map of simplicial presheaves is a \textbf{local inner fibration} (respectively \textbf{local Kan fibration}) if and only if it has the local right lifting property with respect to the inner horn inclusions $\Lambda_{i}^{n} \rightarrow \Delta^{n}, 0 < i < n$ (respectively the horn inclusions  $\Lambda_{i}^{n} \rightarrow \Delta^{n}, 0 \le i \le n$). $\textbf{Local trivial fibrations}$ are defined in a similar manner. If $X$ is a simplicial presheaf so that the map to the terminal sheaf $X \rightarrow *$ has the right lifting property with respect to $\Lambda_{i}^{n} \rightarrow \Delta^{n}, 0 < i < n$, say that $X$ is \textbf{local Joyal fibrant}. Similarly, there is a notion of \textbf{locally Kan fibrant} simplicial presheaves. Note that $X \rightarrow *$ has the sectionwise right lifting property with respect to $\Lambda_{i}^{n} \rightarrow \Delta^{n}, 0 < i < n$ (respectively $\Lambda_{i}^{n} \rightarrow \Delta^{n}, 0 \le i \le n$) if and only if $X$ is a presheaf of quasi-categories (respectively a presheaf of Kan complexes). 

Call a map $f: X \rightarrow Y$ of simplicial presheaves a \textbf{sectionwise Kan fibration} if and only for each $U \in Ob(\mathscr{C})$, $X(U) \rightarrow Y(U)$ is a Kan fibration. There are analagous definitions of \textbf{sectionwise trivial fibrations} and \textbf{sectionwise quasi-fibrations}.
\end{definition}

\begin{lemma}\label{lem2.5} (\cite[Lemma 4.8]{local})
A map of simplicial presheaves $f: X \rightarrow Y$ has the local right lifting property with respect to a finite inclusion of simplicial sets 
$i : K \rightarrow L$ if and only if $\textbf{hom}(L, X) \xrightarrow{(i^{*}, f_{*})} \textbf{hom}(K, X) \times_{\textbf{hom}(K, Y)} \textbf{hom}(L, Y)$ is a local epimorphism. 
\end{lemma}

Throughout the rest of the article fix a Grothendieck site $\mathscr{C}$, and a Boolean localization $ p : s\textbf{Sh}(\mathscr{B}) \rightarrow s\textbf{Sh}(\mathscr{C})$. It is important to note that the Boolean localization is chosen for simplicial sheaves, rather than simplicial presheaves, since a Boolean localization must be a geometric morphism of topoi.

\begin{lemma}\label{lem2.6}
Let $K$ be a finite simplicial set, and $X$ a simplicial presheaf. Then there are natural isomorphisms
\begin{enumerate}
\item{$p^{*}hom(K, L^{2}(X)) \cong hom(K, p^{*}L^{2}(X))$}
\item{$p^{*}\textbf{hom}(K, L^{2}(X)) \cong \textbf{hom}(K, p^{*}L^{2}(X))$}
\item{$L^{2}hom(K, X) \cong hom(K , L^{2}(X))$}
\item{$L^{2}\textbf{hom}(K, X) \cong \textbf{hom}(K, L^{2}(X))$}
\end{enumerate} 
\end{lemma} 

\begin{proof}
1 and 3 are immediate from the fact that $p^{*}, L^{2}$ preserves finite limits and a simplicial set is a colimit of its non-degenerate simplices. The implications 1 $\implies$ 2, 3 $\implies$ 4 are obvious.
\end{proof}

\begin{lemma}\label{lem2.7}
Let $f: X \rightarrow Y$ be a map of simplicial sheaves on a Boolean algebra. Then $f$ has the local right lifting property with respect to inclusion $i : L \rightarrow K$ of finite simplicial sets if and only if it has the sectionwise lifting property with respect to $i$.
\end{lemma}

\begin{proof}
Follows from the axiom of choice for $s\textbf{Sh}(\mathscr{B})$ (\cite[Lemma 3.30]{local}) and \Cref{lem2.5,lem2.6}.
\end{proof}

\begin{lemma}\label{lem2.8} (\cite[Lemma 4.11 and Corollary 4.12.2]{local})
The functors $p^{*}, L^{2}$ both reflect and preserve the property of having the local right lifting property with respect to an inclusion of finite simplicial sets. 
\end{lemma}

\begin{proof}

Follows from \Cref{lem2.5,lem2.6}.
\end{proof}

\begin{definition}\label{def2.9}
A map $f$ of simplicial presheaves is a \textbf{local weak equivalence} if and only if $L^{2}Ex^{\infty}p^{*}L^{2}(f)$ is a sectionwise weak equivalence.
\end{definition}

The intuition behind Boolean localization is that it can be regarded as giving a 'fat' point for a site (for more details see \cite[Section 1]{Boolean}). Thus the definition of local weak equivalence above generalizes the idea of stalkwise weak equivalence in the case of a topos with enough points. This definition of weak equivalence is independent of the choice of Boolean localization. 

\begin{remark}\label{rmk2.10}
 It is clear from the definition of local weak equivalence that $X \rightarrow L^{2}(X)$ is a local weak equivalence. The fact that weak equivalence is independent of the choice of Boolean localization means that if $\mathscr{C}$ is a Boolean site, the choice of Boolean localization can be taken to be the identity. It follows that $p^{*}$ preserves and reflects local weak equivalences.  
\end{remark}

\begin{lemma}\label{lem2.11} (\cite[Corollary 4.28]{local}; also \cite[Corollary 10.9]{Rezk-Boolean}) Let $f : X \rightarrow Y$ be a map of presheaves of Kan complexes. Then $f$ is a local weak equivalence if and only if it is a sectionwise weak equivalence.  
\end{lemma}

\begin{definition}\label{def2.12}
A map of simplicial presheaves $f: X \rightarrow Y$ is said to be a \textbf{local equivalence of fundamental groupoids} if and only if $B\pi(f)$ is a local weak equivalence. There is an analagous notion of \textbf{sectionwise weak equivalences of fundamental groupoids}.
\end{definition}

\begin{lemma}\label{lem2.13}
A map $f: X \rightarrow Y$ of simplicial presheaves of Kan complexes is a local equivalence of fundamental groupoids if and only if $B\pi p^{*}L^{2}(f)$ is a sectionwise weak equivalence.
\end{lemma}

\begin{proof}
By \Cref{lem2.11}, $p^{*}L^{2}(B \pi(f))$ is a sectionwise weak equivalence if and only if $B\pi(f)$ is a local weak equivalence. The result follows from the sequence of natural equivalences
$$
p^{*}L^{2}B \pi(X) \cong p^{*}B L^{2} \pi(X) \simeq p^{*}B L^{2}\pi(L^{2}X) \cong BL^{2}\pi p^{*}L^{2}X.
$$
\end{proof}

\begin{remark}\label{rmk2.14}
Note that if $f: X \rightarrow Y$ is map of Kan complexes then $B\pi(f)$ is a weak equivalence if and only if $\pi(X) \rightarrow \pi(Y)$ is an equivalence of categories. Thus, \Cref{lem2.13} implies that a map $f$ of presheaves of Kan complexes is a local weak equivalence of fundamental groupoids if and only if $\pi(p^{*}L^{2}(X))(b) \rightarrow \pi(p^{*}L^{2}(Y))(b)$ is an equivalence of categories for each $b \in \mathscr{B}$. 
\end{remark}

\begin{definition}\label{def2.15}
Let $$\mathcal{L} : s\textbf{Pre}(\mathscr{C}) \rightarrow s\textbf{Pre}(\mathscr{C})$$ be the functor which applies the usual fibrant replacement functor (i.e. constructed via the small object argument with respect to inner horn inclusions) for the Joyal model structure sectionwise to a simplicial presheaf. If $s\textbf{Pre}(\mathscr{C})^{quasi}, s\textbf{Pre}(\mathscr{C})^{Kan}$ are the full subcategories of $s\textbf{Pre}(\mathscr{C})$ consisting of presheaves of quasi-category and presheaves of Kan complexes, respectively, then sectionwise application of $J$ (as in \Cref{def1.9}) defines a functor $J : s\textbf{Pre}(\mathscr{C})^{quasi} \rightarrow s\textbf{Pre}(\mathscr{C})^{Kan}$. 
\end{definition}

\begin{definition}\label{def2.16}
For a simplicial set $X$, the cardinality of $X$ is defined to be $|X| = \underset{n \in \mathbb{N}}{sup}(|X_{n}|)$. For each simplicial presheaf $X$, and infinite cardinal $\alpha$, say that $X$ is $\alpha$-bounded if $$\underset{U \in Ob(\mathscr{E})}{sup}(|X(U)|) < \alpha$$ Say that a monomorphism $A \rightarrow B$ is $\alpha$-bounded if $B$ is $\alpha$-bounded. 
\end{definition}

\begin{lemma}\label{lem2.17}
There exists an uncountable cardinal $\beth > |Mor(\mathscr{C})|$, so that the following are true: 
\begin{enumerate}
\item{$\mathcal{L}$ preserves filtered colimits. 
}
\item{$\mathcal{L}$ preserves cofibrations.}
\item{Suppose that $\gamma$ is a cardinal so that $\gamma > \beth$. For a simplicial presheaf X, let $\mathcal{F}_{\gamma}(X)$ denote the filtered system of subobjects of X which has cardinality $< \gamma$. The natural map
$$\underset{Y \in \mathcal{F}_{\gamma}(X)}{\underset{\longrightarrow}{\lim}}\mathcal{L}(Y) \rightarrow \mathcal{L}(X)$$ 
is an isomorphism.
}
\item{if $|X| \leq 2^{\lambda}$, where $\lambda \geq \beth$, then $|\mathcal{L}(X)| \leq 2^{\lambda}$.}
\item{$\mathcal{L}$ preserves pullbacks.}
\end{enumerate}
\end{lemma}

\begin{proof}
By arguing sectionwise, this is the same argument as \cite[Theorem 4.8]{J1}. 
\end{proof}

\section{Existence of the Model Structure}\label{sec3}

\begin{definition}\label{def3.1}
Define a map $f: X \rightarrow Y$ of simplicial presheaves on $\mathscr{C}$ to be a \textbf{sectionwise Joyal equivalence} if and only if $X(U) \rightarrow Y(U)$ is a Joyal equivalence for each $U \in Ob(\mathscr{C})$. Define $f$ to be a \textbf{local Joyal equivalence} if and only if $L^{2}\mathcal{L}p^{*}L^{2}(f)$ is a sectionwise Joyal  equivalence. Note that local Joyal equivalences automatically satisfy the 2 out of 3 property. A \textbf{quasi-injective fibration} is a map that has the right lifting property with respect to maps which are both monomorphisms and local Joyal equivalences.
\end{definition}

\begin{corollary}\label{cor3.2}
The map $X \rightarrow L^{2}(X)$ is a local Joyal equivalence.
\end{corollary}

The following theorem is the main theorem of this paper; the remainder of \Cref{sec3} is devoted to its proof.

\begin{theorem}\label{thm3.3}
There exists a left proper model structure on $s\textbf{Pre}(\mathscr{C})$, so that the weak equivalences  are the local Joyal equivalences, the cofibrations are monomorphisms, and the fibrations are the quasi-injective fibrations. 
\end{theorem}

\begin{lemma}\label{lem3.4}
 Let $\Gamma^{*} : sSet \rightarrow s\textbf{Sh}(\mathscr{C})$ be the composite of the constant simplicial presheaf functor and sheafification. The functors$$p^{*}(- \times \Gamma^{*}(C)),  p^{*}(-) \times \Gamma^{*}(C) : s\textbf{Sh}(\mathscr{C}) \rightarrow s\textbf{Sh}(\mathscr{B})$$ are naturally isomorphic for arbitrary simplicial set C.    
\end{lemma}

\begin{proof}
Follows easily by adjunction.
\end{proof}

\begin{lemma}\label{lem3.5}
There is a natural isomorphism $p^{*}L^{2}\mathcal{L} \cong L^{2}\mathcal{L}p^{*}L^{2}$. In particular, $f$ is a local Joyal equivalence if and only if $p^{*}L^{2}\mathcal{L}(f)$ is a sectionwise Joyal equivalence.
\end{lemma}

\begin{proof}
Since $p^{*}, L^{2}$ commute with colimits, by the construction of $\mathcal{L}$, it suffices to show that $p^{*}L^{2}E_{1}(X) \cong L^{2}E_{1}p^{*}L^{2}(X)$ naturally, where $E_{1}$ the pushout of presheaves: 
$$
\xymatrix{
\coprod_{\Lambda_{k}^{n} \subset \Delta^{n}} (hom(\Lambda_{k}^{n}, X) \times \Lambda_{k}^{n}) \ar[r]^>>>>>>{ev} \ar[d] & X \ar[d] \\
\coprod_{\Lambda_{k}^{n} \subset \Delta^{n}} (hom(\Lambda_{k}^{n}, X) \times \Delta^{n}) \ar[r] & E_{1}
}
$$
where the coproducts are indexed over the set of all inner horn inclusions $\Lambda_{k}^{n} \subset \Delta^{n}$, and $ev$ is the evaluation map. Thus, by \Cref{lem2.6,lem3.4} and the fact that sheafification commutes with finite limits, $p^{*}L^{2}(E_{1})$ is naturally isomorphic to the sheaf pushout
$$
\xymatrix{
\coprod_{\Lambda_{k}^{n} \subset \Delta^{n}} (hom(\Lambda_{k}^{n}, p^{*}L^{2}X) \times \Gamma^{*}(\Lambda_{k}^{n})) \ar[r]^>>>>>{ev} \ar[d] & p^{*}L^{2}X \ar[d] \\
\coprod_{\Lambda_{k}^{n} \subset \Delta^{n}} (hom(C, p^{*}L^{2}X) \times \Gamma^{*}(\Delta^{n})) \ar[r] & S
}
$$ 
which is  naturally isomorphic to $L^{2}E_{1}p^{*}L^{2}(X)$, as required. 
\end{proof}

\begin{lemma}\label{lem3.6}Let $X$ be a presheaf of quasi-categories.

\begin{enumerate}
\item{If $X$ is a sheaf of quasi-categories on $\mathscr{B}$, then the natural map $J(X) \rightarrow L^{2}J(X)$ is an isomorphism.}
 \item{For $n \in \mathbb{N}$, let $E_{n}$ denote the set of edges $\Delta^{1} \rightarrow \Delta^{n}$. For each $e \in E_{n}$
form the pullback: 
$$ \xymatrix{
P_{n}^{e} \ar[d] \ar[r]^{\phi_{e}} 
& hom(\Delta^{n}, X)  \ar[d]^{e^{*}} \\
hom(sk_{2}(I), X)  \ar[r]_{i^{*}} & hom(\Delta^{1}, X)}
$$
where $e^{*}, i^{*}$ are induced by inclusion. 
The n-simplices of $J(X)$ are equal to the presheaf-theoretic image of 
$\bigcap_{e \in E_{n}} (P_{n}^{e}) \xrightarrow{\phi} hom(\Delta^{n}, X)$ induced by the $\phi_{e}$'s.}
\item{The $J$ functor commutes with $p^{*}L^{2}$ for presheaves of quasi-categories.}
\end{enumerate}

\end{lemma}

\begin{proof}

First, suppose that $X$ is a sheaf on a Boolean site. Then
$L^{2}J(X)$ is a locally Kan simplicial presheaf, and hence is sectionwise Kan by \Cref{lem2.7}. Furthermore, sheafification preserves injections, so there is a diagram:
$$
\xymatrix{
J(X) \ar@{>->}[r] \ar[d] & X \ar@{=}[d] \\
L^{2}J(X) \ar@{>->}[r] & L^{2}(X) \\
}
$$
Thus $J(X) \rightarrow L^{2}J(X)$ is an inclusion of sub-presheaves of X. But $J(X)$ is the maximal sectionwise Kan subcomplex of X, so that $J(X) = L^{2}(J(X))$

Statement 2 follows immediately from \Cref{lem1.7} and \Cref{thm1.8}.

For the final statement, it is clear that the $P_{n}^{e}$'s are preserved under $p^{*}L^{2}$, since this composite preserves finite limits. Thus, $p^{*}L^{2}(J(X))_{n}$ is isomorphic to the sheaf theoretic image of 
$$\bigcap_{e \in E_{n}} P_{n}^{e} \xrightarrow{\phi} hom(\Delta^{n}, p^{*}L^{2}(X))$$
which is $L^{2}J(p^{*}L^{2}(X))_{n}$ by 2 above. But there is an isomorphism  $L^{2}J(p^{*}L^{2}(X))_{n} \cong Jp^{*}L^{2}(X)_{n}$ by 1 above.

\end{proof}

\begin{theorem}\label{thm3.7}
A map $X \xrightarrow{\phi} Y$ is a local Joyal equivalence if and only if the  map $$J\textbf{hom}(K,  \mathcal{L}(X)) \rightarrow J\textbf{hom}(K, \mathcal{L}(Y))$$ induces a local equivalence of fundamental groupoids for each finite simplicial set $K$. 
\end{theorem}

\begin{proof} \Cref{lem2.6,lem3.6} imply that for each $K$ as above there are isomorphisms: 

$$
p^{*}L^{2}(J\textbf{hom}(K, \mathcal{L}(X))) \cong Jp^{*}L^{2}(\textbf{hom}(K, \mathcal{L}(X))) \cong J\textbf{hom}(K, p^{*}L^{2}\mathcal{L}(X))
$$
so that by \Cref{rmk2.14}, the assertion that the map $\phi$ is a local equivalence of fundamental groupoids is equivalent to 
$$
\pi(J\textbf{hom}(K, p^{*}L^{2}\mathcal{L}(X))(b)) \rightarrow \pi(J\textbf{hom}(K, p^{*}L^{2}\mathcal{L}(Y))(b))
$$
being an equivalence of groupoids for all $b \in Ob(\mathscr{B})$. But since $p^{*}L^{2}\mathcal{L}(X)$ and $p^{*}L^{2}\mathcal{L}(Y)$ are sheaves of quasi-categories, this condition is equivalent to $p^{*}L^{2}\mathcal{L}(\phi)$ being a sectionwise Joyal equivalence, by \Cref{lem1.10}. 
\end{proof}

\begin{remark}\label{rmk3.8}
The preceding theorem shows that the definition of local Joyal equivalence is independent of the Boolean localization chosen. To see this, note that local equivalence of fundamental groupoids can be phrased in terms of isomorphisms of sheaves of homotopy groups (c.f. \cite[pg.63-64]{local}). In particular the Boolean localization can be taken to be the identity if $\mathscr{C}$ is a Boolean site. Thus, \Cref{lem3.4} implies that $p^{*}$ and $L^{2}$ both preserve and reflect local Joyal equivalences. 
\end{remark}

\begin{corollary}\label{cor3.9}
A sectionwise Joyal equivalence $X \rightarrow Y$ of simplicial presheaves is a local Joyal equivalence.
\end{corollary}

\begin{corollary}\label{cor3.10}
A local Joyal equivalence between sheaves of quasi-categories on a Boolean site is a sectionwise Joyal equivalence. 
\end{corollary}

\begin{corollary}\label{cor3.11} A map $f: X \rightarrow Y$ of presheaves of quasi-categories is a local Joyal equivalence if and only if $p^{*}L^{2}(f)$ is a sectionwise Joyal equivalence.
\end{corollary}

\begin{proof}
Suppose that $f$ is a local Joyal equivalence. The map $p^{*}L^{2}(X) \rightarrow \mathcal{L}p^{*}L^{2}(X) \rightarrow L^{2}\mathcal{L}p^{*}L^{2}(X)$ is a local Joyal equivalence in $s\textbf{Pre}(\mathscr{B})$ by \Cref{cor3.2,cor3.9}. Furthermore, $L^{2}\mathcal{L}p^{*}L^{2}(f)$ is a sectionwise, and hence local Joyal equivalence. Thus, the commutative diagram

$$
\xymatrix{
p^{*}L^{2}(X) \ar[r] \ar[d] & L^{2}\mathcal{L} p^{*}L^{2}(X) \ar[d]^{L^{2}\mathcal{L}p^{*}L^{2}(f)} \\
p^{*}L^{2}(Y) \ar[r] & L^{2} \mathcal{L}p^{*} L^{2}(Y) \\
}
$$
and the 2 out of 3 property imply that $p^{*}L^{2}(f)$ is a local Joyal equivalence in s$\textbf{Pre}(\mathscr{B})$. But a local Joyal equivalence between sheaves of quasi-categories on $\mathscr{B}$ is a sectionwise Joyal equivalence. 

The converse is similar, but easier.
\end{proof}

\begin{corollary}\label{cor3.12}
 A map $f$ of simplicial presheaves is a local Joyal equivalence if and only if $\mathcal{L}(f)$ is a local Joyal equivalence. 
\end{corollary}

\begin{corollary}\label{cor3.13}
Let $X_{\alpha} \rightarrow Y_{\alpha}$ be natural transformations consisting of local Joyal equivalences of presheaves of quasi-categories which is indexed by some filtered category $J$. Then the induced map 
$$\underset{\alpha \in J}{\underset{\longrightarrow}{\lim}} X_{\alpha} \rightarrow \underset{\alpha \in J}{\underset{\longrightarrow}{\lim}} Y_{\alpha}$$
is a local Joyal equivalence. 
\end{corollary}

\begin{lemma}\label{lem3.14}
A map $f: X \rightarrow Y$ of quasi-categories is a quasi-fibration if and only if it is an inner fibration and there exists a lift in each diagram of the form
$$
\xymatrix{
\Delta^{0} \ar[r] \ar[d]_{d^{0}} & J(X) \ar[d] \\
\Delta^{1} \ar[r] \ar@{.>}[ur] & J(Y) 
}
$$

Furthermore, $p^{*}L^{2}$ both preserves and reflects the property of being a sectionwise quasi-fibration of presheaves of quasi-categories.

\end{lemma}

\begin{proof}
The first statement is immediate from \cite[Corollary 2.4.6.5]{Lurie}; see also \cite[Definition 4.2]{Joyal-quasi-cat}. The second statement follows from the first and \Cref{lem2.7,lem2.8,lem3.6}
\end{proof}

\begin{lemma}\label{lem3.15}
A local trivial fibration $f: X \rightarrow Y$ is a local Joyal equivalence. Suppose that $f: X \rightarrow Y$ is a sectionwise quasi-fibration of presheaves of quasi-categories. Then $f$ is a local Joyal equivalence if and only if it is a locally trivial fibration.
\end{lemma}

\begin{proof}
First note that if $f$ is any local trivial fibration, then $p^{*}L^{2}(f)$ is a sectionwise trivial fibration by \Cref{lem2.7,lem2.8}. Thus it is a sectionwise and hence local Joyal equivalence. 

Now, suppose that $f$ is a local Joyal equivalence of presheaves of quasi-categories and a sectionwise quasi-fibration. By \Cref{lem3.14}, $p^{*}L^{2}(X) \rightarrow p^{*}L^{2}(Y)$ is a sectionwise quasi-fibration. By \Cref{cor3.11}, it is also sectionwise Joyal equivalence so  $p^{*}L^{2}(X) \rightarrow p^{*}L^{2}(Y)$ is a sectionwise, and hence locally, trivial  fibration. The result follows from \Cref{lem2.8}.
\end{proof}

\begin{example}\label{exam3.16}
This example gives the construction of the quasi-fibration replacement for a map $f: X \rightarrow Y$ of presheaves of quasi-categories. 
Form the diagram of simplicial presheaves 
$$
\xymatrix{
X \times_{Y} \textbf{hom}(I, Y) \ar[d]_{d_{0*}} \ar[r]^>>>>>{f_{*}} & \textbf{hom}(I, Y) \ar[r]^>>>>{d_{1}} \ar[d]^{d_{0}} & Y \\
X \ar[r]_{f} \ar[ur]^{sf}  & Y & 
 }
$$
Since $\textbf{hom}(I, Y)$ is a path object for the Joyal model structure and $Y$ is a presheaf of quasi-categories, and $d_{0}$ is a sectionwise trivial fibration, so that $d_{0*}$ is a sectionwise trivial fibration. The section $s$ of $d_{0}$ induces a section $s_{*}$ of $d_{0*}$, and 
$$
(d_{1}f_{*})s_{*} = d_{1}sf = f
$$ 
Finally, there is a pullback diagram of presheaves:
$$
\xymatrix{
X \times_{Y} \textbf{hom}(I, Y)\ar[r]^{f_{*}} \ar[d]_{(d_{0*}, d_{1}f_{*})} & \textbf{hom}(I, Y) \ar[d]^{(d_{0}, d_{1})} \\
X \times Y \ar[r]_{f \times 1} & Y \times Y  
}
$$
and the projection map $pr_{R}: X \times Y \rightarrow Y$ is a sectionwise quasi-fibration, since $X$ is a presheaf of quasi-categories. Thus, $pr_{R}(d_{0*}, d_{1}f_{*}) = d_{1}f_{*}$ is a sectionwise quasi-fibration. Write $Z_{f} = X \times_{Y} Y^{I}$, and $\pi = d_{1}f_{*}$. Then $\pi$ is a functorial replacement of $f$ by a quasi-fibration, and there is a diagram
$$
\xymatrix{
X \ar[r]^{s_{*}} \ar[dr]_{f} & Z_{f} \ar[d]_{\pi} \ar[r]^{(d_{0})_{*}} & X \\
& Y & }
$$
where $(d_{0})_{*}$ is a trivial fibration and $(d_{0})_{*} \circ s_{*} = id_{X}$.
\end{example}

\begin{remark}\label{rmk3.17}
An analagous construction to that above produces the sectionwise Kan fibration replacement of a  map of presheaves of Kan complexes. Taking pullbacks gives a functorial Kan fibration replacement for all simplicial presheaf maps. However, this technique does not work for the local Joyal model structure, since the Joyal model structure is not right proper. 
\end{remark}

\begin{lemma}\label{lem3.18}

Let $\alpha$ be a regular cardinal so that $\alpha > |Mor(\mathscr{E})|$ and let $C \subseteq A$ be an inclusion of simplicial presheaves, so that $C$ is $\alpha$-bounded and $A$ is a presheaf of quasi-categories. Then there exists an $\alpha$-bounded presheaf of quasi-categories $B$ so that $C \subseteq B \subseteq A$.
\end{lemma}

\begin{proof}
The set of lifting problems 
$$
\xymatrix{
\Lambda_{i}^{n} \ar[r] \ar[d] & C(U) \\
\Delta^{n} \ar@{.>}[ur]& \\ }
$$
for $U \in Ob(\mathscr{E})$ is $\alpha$-bounded and can be solved over $A$. Furthermore, since $A$ is a colimit of its $\alpha$-bounded subobjects there is a subobject $B_{1}$ of  $A$ so that $C \subseteq B_{1}$, all lifting problems as above can be solved over $B_{1}$, and $B_{1}$ is $\alpha$-bounded. Repeating this procedure countably many times produces an ascending sequence 
$$
B_{1} \subseteq B_{2} \subseteq \cdots \subseteq B_{n} \subseteq \cdots
$$ 
Set $B = \cup_{i=1}^{\infty} B_{i}$.
\end{proof}

\begin{lemma}\label{lem3.19} (see \cite[Theorem 5.2]{local})
Suppose that $\alpha$ is a regular cardinal, so that $\alpha > |Mor(\mathscr{C})|$. Suppose that there is a diagram of monomorphisms of simplicial presheaves of quasi-categories as below: 
$$
\xymatrix{
 & X \ar[d] \\
A \ar[r] & Y
}
$$
where $A$ is $\alpha$-bounded, and $X \rightarrow Y$ is a local Joyal equivalence. Then there exists an $\alpha$-bounded presheaf of quasi-categories $B$, so that $A \subseteq B \subseteq Y$ and $B \cap X \rightarrow B$ is a local Joyal equivalence. 
\end{lemma}

\begin{proof}
If $B$ is a presheaf of quasi-categories, write $\pi_{B} : Z_{B} \rightarrow B$ for the natural quasi-fibration replacement of $B \cap X \rightarrow B$ (as explained in \Cref{exam3.16}). By \Cref{lem3.15} $B \cap X \rightarrow B$ is a local Joyal equivalence equivalence if and only if $\pi_{B}$ is a local trivial fibration. Now, suppose there is a lifting problem 

$$
\xymatrix{
\partial \Delta^{n} \ar[r] \ar[d] & Z_{A}(U) \ar[d] \\
\Delta^{n} \ar[r] \ar@{.>}[ur] & A(U) \\
}
$$
Then this lifting problem can be solved locally over some covering $\{ U_{i} \rightarrow U\}$ having at most $\alpha$ elements. There is an identification 
$$
\underset{|B| < \alpha}{\underset{\longrightarrow}{\lim}}Z_{B}= Z_{Y}
$$ 
Thus, it follows from the regularity assumption on $\alpha$ there is an $\alpha$-bounded $A' \subseteq Y$, $A \subseteq A'$, that can solve the lifting problem above. The set of all such lifting problems is $\alpha$-bounded. Thus, there is a $\alpha$-bounded presheaf of quasi-categories $B_{1} \subseteq Y$ such that each lifting problem as above can be solved over $B_{1}$ by \Cref{lem3.18}. Repeating this procedure countably many times produces an ascending sequence of presheaves of quasi-categories
$$B_{1} \subseteq B_{2} \cdots \subseteq B_{n} \cdots$$
such that all lifting problems:
$$
\xymatrix{
\partial \Delta^{n} \ar[r] \ar[d] & Z_{B_{i}}(U) \ar[d] \\
\Delta^{n} \ar[r] \ar@{.>}[ur] & B_{i}(U) \\
}
$$
can be solved locally over $B_{i+1}$. Put $B = \cup B_{i}$. Then $B$ is $\alpha$-bounded by the regularity of $\alpha$. Furthermore $B$ is a presheaf of quasi-categories. Since the construction of $Z_{B}$ commutes with filtered colimits, $Z_{B} \rightarrow B$ is a local trivial fibration, as required.  
\end{proof}

\begin{lemma}\label{lem3.20}
Let $\beta > 2^{\beth}$, where $\beth$ is as in \Cref{lem2.17}. Also, put $\alpha = 2^{\beta} + 1$, so that $\alpha$ is a regular cardinal since it is a successor. Suppose that there is a diagram of monomorphisms of simplicial presheaves as below:
$$
\xymatrix{
 & X \ar[d] \\
A \ar[r] & Y
}
$$
where $A$ is $\alpha$-bounded, and $X \rightarrow Y$ is a local Joyal equivalence. Then there exist $\alpha$-bounded simplicial presheaves $A', B'$, so that 
\begin{enumerate}
\item{ $\mathcal{L}(A) \subseteq A' \subseteq \mathcal{L}(Y)$, $\mathcal{L}(X) \cap A' \rightarrow A'$ is a local Joyal equivalence.}
\item{$A \subseteq B' $, $A' \subseteq \mathcal{L}(B')$}
\end{enumerate}

\end{lemma}

\begin{proof}
Since $\mathcal{L}$ preserves monomorphisms, \Cref{lem2.17} implies that there is a diagram of $\alpha$-bounded monomorphisms:
$$
\xymatrix{
 & \mathcal{L}(X) \ar[d] \\ 
\mathcal{L}(A) \ar[r] & \mathcal{L}(Y)
}
$$
Hence there is an $A'$ with the desired properties by \Cref{lem3.19}. Now, note that by \Cref{lem2.17}: 
$$\underset{Y \in \mathcal{F}_{\alpha}(X)}{\underset{\longrightarrow}{\lim}}\mathcal{L}(Y) \cong \mathcal{L}(X)$$ 
Furthermore, the set of elements: $$\{(x, U): x \in A'(U) - \mathcal{L}(A)(U) , U \in Ob(\mathscr{E})\}$$ is $\alpha$-bounded, so there exists an $\alpha$-bounded object $B'$ with the desired properties. 
\end{proof}

\begin{theorem}\label{thm3.21}
Let $\alpha$ be as in \Cref{lem3.20}. Suppose that there is diagram of monomorphisms of simplicial presheaves of  as below:
$$
\xymatrix{
 & X \ar[d] \\
A \ar[r] & Y
}
$$
where $A$ is $\alpha$-bounded, and $X \rightarrow Y$ is a local Joyal equivalence. Then there exists a $\alpha$-bounded subobject $B, A \subseteq B \subseteq Y$, so that $B \cap X \rightarrow B$ is a local Joyal equivalence. 
\end{theorem}

\begin{proof}
For each $n \in \mathbb{N}$, define $\alpha$-bounded objects $A_{n}, B_{n}$ inductively, so that the following properties hold: 
\begin{enumerate}
\item{$\mathcal{L}(B_{n'}) \subseteq A_{n} \subseteq \mathcal{L}(Y)$ for all $n' < n$, $\mathcal{L}(X) \cap A_{n} \rightarrow A_{n}$ is a local Joyal equivalence}
\item{$A \subseteq B_{n} \subseteq Y \, \, A_{n} \subseteq \mathcal{L}(B_{n})$}
\end{enumerate}
Start the induction by setting $A_{0} = B_{0} = A$. In general, having defined $A_{n'}, B_{n'}$ for $n' < n$ apply \Cref{lem3.20} to the diagram
$$
\xymatrix{
& \ar[d] X \\
B_{n-1} \ar[r] & Y }
$$

to produce $A_{n}, B_{n}$. 

Let $B = \underset{n \in \mathbb{N}}{\underset{\longrightarrow}{\lim}}B_{n}$. 
B is $\alpha$-bounded by the regularity of $\alpha$. Now, note that by \Cref{lem2.17}, for $X'$ a subobject of Y, there are natural isomorphisms:

$$
\mathcal{L}(B \cap X') \cong \underset{n \in \mathbb{N}}{\underset{\longrightarrow}{\lim}}\mathcal{L}(B_{n} \cap X') 
\cong \underset{n \in \mathbb{N}}{\underset{\longrightarrow}{\lim}}\mathcal{L}(B_{n}) \cap \mathcal{L}(X') 
\cong \underset{n \in \mathbb{N}}{\underset{\longrightarrow}{\lim}}(A_{n} \cap \mathcal{L} (X')) 
$$
so that $\mathcal{L}(B \cap X) \rightarrow \mathcal{L}(B)$ is a local Joyal equivalence by \Cref{cor3.13}. Thus, the map $B \cap X \rightarrow X$ is a local Joyal equivalence by \Cref{cor3.12}, as required. 
\end{proof}

\begin{lemma}\label{lem3.22}
Let $\alpha$ be a cardinal as in \Cref{lem3.20} and \Cref{thm3.21}. There a map $f$ has the right lifting property with respect to all maps which are cofibrations (respectively local Joyal equivalences and cofibrations) if and only it has the right lifting property with respect to all $\alpha$-bounded cofibrations (respectively $\alpha$-bounded local Joyal equivalences and cofibrations).
\end{lemma}

\begin{proof}
For cofibrations, this is just \cite[Theorem 5.6]{local}. For cofibrations that are local Joyal equivalences, use Theorem 3.21 and the method of \cite[Lemma 5.4]{local}.
\end{proof}

\begin{lemma}\label{lem3.23}
A map $f: X \rightarrow Y$ of simplicial presheaves which has the right lifting property with respect to all cofibrations is a local Joyal equivalence and a quasi-injective-fibration. 
\end{lemma}

\begin{proof}
The map $f$ is a quasi-injective fibration by definition. $f$ is also a sectionwise trivial Kan fibration, and hence a local trivial Kan fibration. Conclude using \Cref{lem3.15}.
\end{proof}

\begin{lemma}\label{lem3.24}
Consider a pushout diagram of simplicial presheaves

$$
\xymatrix{
A \ar[r]^{\alpha} \ar[d]_{\beta} & B \ar[d]^{\beta'} \\
C \ar[r]_{\alpha'} & B \cup_{A} C  
}
$$
where $\alpha$ is a cofibration. Then $\beta'$ is a local Joyal equivalence if $\beta$ is. 

\end{lemma}

\begin{proof}
In the case that $A. B$ and $C$ are sheaves of quasi-categories on the Boolean algebra $\mathscr{B}$, this is immediate from the left properness of the Joyal model structure, \Cref{rmk3.8} and \Cref{cor3.11}.   

Now, suppose A, B, C, D are arbitrary simplicial presheaves. In the diagram below each of the vertical maps are sectionwise Joyal equivalences: 
$$
\xymatrix{
B \ar[d] & A \ar[l] \ar[r] \ar[d] &  C \ar[d] \\
\mathcal{L}(B) & \mathcal{L}(A) \ar[l] \ar[r] & \mathcal{L}(C) \\
}
$$

The gluing lemma (\cite[Lemma 1.8.8]{GJ2}) implies that the induced map $B \cup_{A} C \rightarrow \mathcal{L}(B) \cup_{ \mathcal{L}(A)} \mathcal{L}(C)$ is a sectionwise and hence local Joyal equivalence. Thus the diagram 
$$
\xymatrix{
A \ar[r]^{\alpha} \ar[d]_{\beta} & B \ar[d]^{\beta'} \\
C \ar[r]_{\alpha'} & B \cup_{A} C  
}
$$  
is local Joyal equivalent to:
$$
\xymatrix{
\mathcal{L}(A) \ar[r]^{\mathcal{L}(\alpha)} \ar[d]_{\mathcal{L}(\beta)} & \mathcal{L}(B) \ar[d]^{s} \\
\mathcal{L}(C) \ar[r]_>>>>>{s'} & \mathcal{L}(B) \cup_{\mathcal{L}(A)} \mathcal{L}(C) \\
}
$$  
Since $p^{*}L^{2}$ preserves pushouts and cofibrations, the case of sheaves of quasi-categories on $\mathscr{B}$ implies that $p^{*}L^{2}(s)$ is a local Joyal equivalence. Thus, so is $s'$, since local Joyal equivalences are reflected by Boolean localization (\Cref{rmk3.8}). 

\end{proof}

\begin{lemma}\label{lem3.25}
Let $f: X \rightarrow Y$ be a map of simplicial presheaves. Then it can be factored as 
$$
\xymatrix{
& Z \ar[dr]^{p} & \\
X \ar[ur]^{i} \ar[dr]_{q} \ar[rr]^{f} && Y \\
& W \ar[ur]_{j} & }
$$ 
where 
\begin{enumerate}
\item{$i$ is a local Joyal equivalence and a cofibration and $p$ is a quasi-injective fibration.}
\item{$j$ is a cofibration and $q$ is a quasi-injective fibration and local Joyal equivalence.}
\end{enumerate}

\end{lemma}

\begin{proof}
For the first factorization choose a cardinal $\lambda > 2^{\alpha}$, where $\alpha$ is chosen as in \Cref{lem3.22} to solve all lifting problems
$$
\xymatrix{
A \ar[r] \ar[d]_{i} & X \ar[d]^{f} \\
B \ar[r] \ar@{.>}[ur] & Y
}
$$ 
where $i$ is an $\alpha$-bounded trivial cofibration. The result follows from the fact that local Joyal equivalences are preserved under pushout, which is \Cref{lem3.24}. 

The second statement is proven in a similar manner, using \Cref{lem3.23}. 

\end{proof}

\begin{proof}[Proof of Theorem 3.3]

CM1, CM2 are trivial to verify. CM3 follows from the definition of local Joyal equivalences. CM5 is \Cref{lem3.25}.

Let $f: X \rightarrow Y$ be a map which is a local Joyal equivalence and a quasi-injective fibration. The next paragraph will show that $f$ has the right lifting property with respect to cofibrations. By \Cref{lem3.25}, the map has a factorization 
$$
\xymatrix{
X \ar[r]^{g} \ar[dr]_{f} & W \ar[d]^{h} \\
 & Y \\
}
$$ 
where $h$ has the right lifting property with respect to all cofibrations, and is therefore a local Joyal equivalence, and g is a cofibration. Hence by the 2 out of 3 property, g is a local Joyal equivalence and a cofibration. Thus there is a lifting in the below diagram
$$
\xymatrix{
X \ar[r]^{id} \ar[d]_{g} & X \ar[d]^{f} \\
W \ar[r]_{h} \ar@{.>}[ur]^{n} & Y
}
$$
Finally, the diagram 
$$
\xymatrix{
X \ar[r]^{g} \ar[d]_{f} & W \ar[d]^{h} \ar[r]^{n} & X \ar[d]^{f} \\
Y \ar[r]_{id} & Y \ar[r]_{id} & Y \\
& & 
}
$$
shows that $f$ is a retract of $h$ and hence $f$ has the right lifting property with respect to all cofibrations (since right lifting property is preserved under retracts), as required. (The argument is standard; for instance, see \cite[Theorem 5.8]{local})

\end{proof}

\section{The Model Structure on Simplicial Sheaves, and an Example}

\begin{theorem}\label{thm4.1}
The category $s\textbf{Sh}(\mathscr{C})$ along with the class of local Joyal equivalences, monomorphisms, and quasi-injective fibration forms a left proper model structure.

Furthermore, there is a Quillen adjunction
$$
L^{2} : s\textbf{Pre}(\mathscr{C}) \leftrightarrows s\textbf{Sh}(\mathscr{C}): i
$$
where $i$ is the inclusion of sheaves into presheaves and $L^{2}$ is sheafification. 
\end{theorem}

\begin{proof}
The associated sheaf functor preserves and reflects local Joyal equivalences, by \Cref{rmk3.8}, and it also preserves cofibrations. Hence, the inclusion functor preserves quasi-injective fibrations. Thus the functors form a Quillen pair. The unit map of the adjunction $X \rightarrow L^{2}(X)$ is a local Joyal equivalence, and the counit map is the identity. Thus, the second statement follows from the first, and it suffices to prove the first statement. 

 Axiom CM1 follows from completeness and cocompleteness of the sheaf category. Axioms CM2-CM4 follow from the corresponding statements for simplicial presheaves. 
 Let $\alpha$ be a cardinal as in \Cref{lem3.22}. Then choose a regular cardinal $\beta$ so that $L^{2}(f)$ is $\beta$ bounded for each $\alpha$-bounded trivial cofibration. Then a map $f$ is a quasi-injective fibration if and only if it has the right lifting property with respect to all $\beta$-bounded trivial cofibration. Doing a small object argument of size $2^{\beta}$ as in Lemma 3.25 gives one half of CM5. The other half has an analagous proof. 

Left properness comes from the corresponding statement for simplicial presheaves, as well as the fact that $X \rightarrow L^{2}(X)$ is a local Joyal equivalence.
\end{proof}

It is asserted in
\cite[Theorem 10.6]{Rezk-Boolean} that for the Jardine model structure on $s\textbf{Sh}(\mathscr{B})$, the injective fibrations are precisely the sectionwise Kan fibrations and the trivial injective fibrations are the sectionwise trivial fibrations. The following analogue of this theorem is true for the local Joyal model structure. 

\begin{theorem}\label{thm4.2}
The local Joyal model structure on $s\textbf{Sh}(\mathscr{B})$ has the following description
\begin{enumerate}
\item{The cofibrations are monomorphisms.}
\item{The fibrations are sectionwise quasi-fibrations.}
\item{The weak equivalences are local Joyal equivalences.} 
\end{enumerate}
\end{theorem}

\begin{lemma}\label{lem4.3}
A map $f$ of $s\textbf{Sh}(\mathscr{B})$ which is a sectionwise quasi-fibration and a local Joyal equivalence is a sectionwise trivial fibration.
\end{lemma}

\begin{proof}
Suppose that that $X \rightarrow Y$ is a map with the properties stated. Let $\mathcal{L}(X) \rightarrow X' \rightarrow \mathcal{L}(Y)$ be the functorial factorization of \Cref{exam3.16} and consider the pullback

$$
\xymatrix{
P \ar[d] \ar[r] & X' \ar[d] \\
Y \ar[r] & \mathcal{L}(Y)
}
$$
The right vertical map is a sectionwise trivial Kan fibration by \Cref{lem2.7,lem3.15}. Hence the left vertical map is as well. Thus by two out of three, all of the maps in the preceding commutative square are sectionwise Joyal equivalences. 

Consider the diagram
$$
\xymatrix{
X \ar[rd]^>>>{\phi} \ar[rrd] \ar[rdd]  & &  \\
& P \ar[r] \ar[d] & X' \ar[d] \\
& Y \ar[r] & \mathcal{L}(Y)
}
$$
where $\phi$ is the induced map. By two out of three $X \rightarrow P$, and hence $X \rightarrow Y$ are sectionwise Joyal equivalence. But then $X \rightarrow Y$ is a sectionwise Joyal equivalence and a sectionwise quasi-fibration, from which the result follows. 
\end{proof}

\begin{lemma}\label{lem4.4}
 Let $\mathcal{D}$ denote a set of generating trivial cofibrations for the Joyal model structure. Denote by $\mathcal{C}$ the set of maps that are retracts of transfinite composites of pushouts of maps of the form $y(b) \times \phi$, where $\phi \in \mathcal{D}$. Here $y$ denotes the yoneda embedding.
\begin{enumerate}
\item{ A map $f: X \rightarrow Y$ of simplicial sheaves on $\mathscr{B}$ admits a factorization as an element of $\mathcal{C}$ followed by a sectionwise quasi-fibration.}
\item{The elements of $\mathcal{C}$ are precisely maps which are trivial cofibrations for the local Joyal model structure on $s\textbf{Sh}(\mathscr{B})$.}
\end{enumerate}
\end{lemma}

\begin{proof}
For the first statement of the lemma, choose a cardinal $\lambda > |\mathscr{B}|$, so that $\lambda > |B| $ for each $\phi : A \rightarrow B$ in $\mathcal{D}$. Do a small object argument of size $2^{\lambda} + 1$ to solve all lifting problems 
$$
\xymatrix{
y(b) \times A \ar[r] \ar[d]_{id_{y(b)} \times \phi} & X \ar[d] \\
y(b) \times B \ar[r] \ar@{.>}[ur] & Y }
$$
where $b \in \mathscr{B}$. This gives the required factorization.

The fact that each member of $\mathcal{C}$ is a local trivial cofibration follows from the left properness of the local Joyal model structure and the fact that local Joyal equivalences are preserved by filtered colimits. 

For the converse, let $f$ be a trivial cofibration for the local Joyal model structure. Factor $f = g \circ h$ where $g$ is a sectionwise quasi-fibration and $h \in \mathcal{C}$. The map $g$ is a sectionwise trivial fibration by \Cref{lem4.3}. Thus
\cite[Lemma 10.14]{Rezk-Boolean} (i.e. CM5 for the model structure of \cite[Theorem 10.6]{Rezk-Boolean}), can be used to show that $f$ is a retract of $h$, so that $f \in \mathcal{C}$.

\end{proof}

\begin{proof}[Proof of \Cref{thm4.2}]
Once it is proven that that the above description gives a model structure on $s\textbf{Sh}(\mathscr{B})$, is immediate that it coincides with the local Joyal model structure.

CM1-CM3 are trivial. The factorization of a map as a trivial cofibration followed by a sectionwise quasi-fibration follows from \Cref{lem4.4} and CM5 for the local Joyal model structure. The factorization of a map as a cofibration followed by a trivial fibration follows from CM5 for the model structure of \cite[Theorem 10.6]{Rezk-Boolean}, and \Cref{lem3.15}. One half of CM4 follows from \Cref{lem4.4}. The other half follows from CM4 for the model structure of \cite[Theorem 10.6]{Rezk-Boolean} and \Cref{lem4.3}.

\end{proof}

\bibliography{database}

\end{document}